\documentclass{article} 
\usepackage{mathmag}

\usepackage{amsmath,amsthm}     
\usepackage{graphicx}     
\usepackage{hyperref} 
\usepackage{url}
\usepackage{amsfonts}

\usepackage{enumitem}
\usepackage{float}
\usepackage{multicol}
\usepackage{tikz}

\theoremstyle{theorem}
\newtheorem{theorem}{Theorem}

\theoremstyle{theorem}
\newtheorem{proposition}{Proposition}

\theoremstyle{theorem}
\newtheorem{question}{Question}

\theoremstyle{definition}

\allowdisplaybreaks

\makeatletter
\@addtoreset{footnote}{page}
\makeatother

\begin{document}

\title{Eulerian 2-Complexes}

\author{Richard H.\ Hammack\\               
{\scriptsize Department of Mathematics and Applied Mathematics}\\
{\scriptsize Virginia Commonwealth University}\\            
{\scriptsize rhammack@vcu.edu}\\
\bigskip
Paul C.\ Kainen\\
{\scriptsize Department of Mathematics and Statistics}\\    
{\scriptsize Georgetown University}\\           
{\scriptsize kainen@georgetown.edu}\\
}                     .

\maketitle

\noindent  
A celebrated theorem of graph theory, known to mathematicians
and non-mathematicians alike, is what we call the {\it Eulerian equivalence}:

\begin{theorem}
\label{Theorem:Euler}
The following are equivalent for a connected graph~$G$.
\begin{itemize}[itemsep=0pt, align=right,topsep=2pt, labelwidth=1.5em,leftmargin=3.5em]
\item [{\em (i)}] Each vertex of $G$ has positive even degree\hfill {\em (in which case we say $G$ is {\bf even})},
\item [{\em (ii)}] $G$ is an edge-disjoint union of cycles,
\item [{\em (iii)}] $G$ has an Euler tour.
\end{itemize}
\end{theorem}

This is traditionally attributed to Euler's 1735 paper~\cite{EULER},
on what is now called the K\"onigsberg bridge problem (although Euler gave
only a partial proof).
Euler asked whether
one could cross each bridge in a city exactly once, returning to the starting point. 
He modeled the problem with a graph, where each land mass is a vertex and each bridge joining two land masses
is an edge between the corresponding vertices.  A successful traversal of the bridges corresponds to
an {Euler tour} in the graph, which we view informally as a closed route through the graph that traverses each edge once.
(See \cite{CHARTRAND} for a more careful definition and for other standard terms not defined here.)

\begin{figure}[H]
\centering
\includegraphics{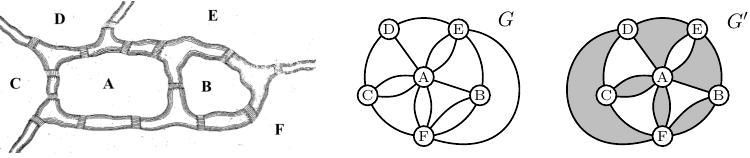}
\caption{Left: A network of bridges. Center: The graph that models the network. Right: A graph that meets the conditions of the Eulerian equivalence. {\small \it arXiv copy, 2023}}
\label{Fig:Bridges}
\end{figure}

For example, consider the map in Figure~\ref{Fig:Bridges} (left),
which is adapted from Euler's paper. This configuration is modeled by
 the graph $G$ (center). This graph has two vertices of odd degree,
so the Eulerian equivalence implies that it has no Euler tour. By contrast, every vertex of the graph $G'$ (right) has
even degree, so $G'$ has an Euler tour. (The reader can easily find one.) Further, $G'$ is the edge-disjoint union of the
six cycles that are the boundaries of the six shaded regions. (The shading is included only to indicate the cycles.
Note that we allow for cycles of length 2.)

This article explores the generalization of the eulerian equivalence from graphs\linebreak (1-complexes) to 2-complexes.
Certain technical issues arise in extending this beyond dimension 2, but the theory
for 2-complexes is, as we will see, very rich.

We assume only a working knowledge of graph theory
and surface topology.

For us, a {\bf 2-complex} is a  triple $K=\big(V(K), E(K), F(K)\big)$ of three finite sets: a {\bf vertex set} $V(K)$ of points; an
{\bf edge set} $E(K)$ whose elements are homeomorphic to the unit interval;
 and a {\bf face set} $F(K)$ of polygons,
each homeomorphic to some closed $n$-gon ($n\geq 3$) in $\mathbb{R}^2$.
In addition, the endpoints of any edge in $E(K)$ are distinct elements of $V(K)$,
and the edges of any polygon in $F(K)$ are distinct edges in $E(K)$.
Further, the elements of $F(K)$ have pairwise non-intersecting
interiors, as do those of $E(K)$. Elements of $V(K)$, $E(K)$ and $F(K)$ are also called
{\bf 0-cells},  {\bf 1-cells} and  {\bf 2-cells}, respectively.
 
The {\bf degree} $\deg(v)$ of a $v\in V(K)$ is the number of edges to which $v$ belongs;
likewise, the degree $\deg(e)$ of $e\in E(K)$ is the number of faces to which $e$ belongs.
A {\bf 1-complex} is a pair $K=(V(K), E(K))$ for which the endpoints of any edge in $E(K)$ are
distinct vertices in $V(K)$. (We regard graphs as 1-complexes.) 
The {\bf 1-skeleton} of  a 2-complex $K=\big(V(K), E(K), F(K)\big)$ is the 1-complex
$K^1:=\big(V(K), E(K)\big)$.
A  {\bf cellular map} $K\to K'$ is a map between complexes 
whose\linebreak restriction to any $k$-cell of $K$ is a homeomorphism onto a $k$-cell  of $K'$
for $k=0,1,2$.

We regard a 2-complex as a topological space in the obvious way,  as
the union of its cells. 
It is a {\bf surface} if each of its points has a neighborhood homeomorphic
to an open disk. Clearly each edge of a surface belongs to exactly two faces.

 \begin{figure}[h]
 \centering
\reflectbox{\includegraphics[width=4in]{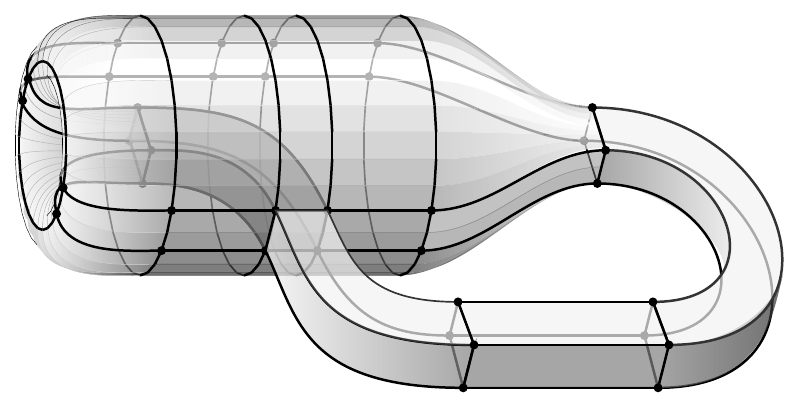}}
 \caption{A 2-complex $K$.}
 \label{Fig:Klein}
 \end{figure}
 
 Figure~\ref{Fig:Klein} shows a 2-complex $K$, embedded in an ambient space $\mathbb{R}^3$.
All of its faces happen to be squares.
 It is not a surface, as some of its edges belong to four faces. We might regard it as the image
 $\varphi(M)$ of a surface $M$ under a cellular map $\varphi$ that is bijective on
$M-M^1$ (so the image of $M$ intersects only along its 1-skeleton).
 
 To illustrate these terms, and also to motivate the main idea of this article, we pose a 
 quick multiple choice question.

 \begin{question} The 2-complex $K$ shown in Figure~\ref{Fig:Klein} is
\begin{itemize}[itemsep=0pt, align=right,topsep=3pt, labelwidth=1.5em,leftmargin=3.5em]
\item [{\em (a)}] the image of a Klein bottle,
\item [{\em (b)}] the image of a torus,
\item [{\em (c)}] the image of a double (two-holed) torus,
\item [{\em (d)}] the image of a triple (three-holed) torus.
\end{itemize}
\end{question}
The question is easy, not because the answer is obvious, but because all four answers are correct.
(There are even four additional correct answers!)
As we will see shortly, each choice is a different ``traversal'' of $K$ by a ``2-dimensional
Euler tour.''

\section{The Euler Equivalence for 2-Complexes}

Let us now adapt conditions (i), (ii), (iii) of the eulerian equivalence to 2-complexes.

The {\it even-degree} condition (i) is replaced with the requirement that
each {\it edge} of a 2-complex have positive even degree. Such a
2-complex is called {\bf even} or {\bf Eulerian}. 

Extending condition (ii) to $2$-complexes requires generalizing the notion of a cycle in a graph.
A {\bf cycle} is an even 1-complex that is not an
edge-disjoint union of two even 1-complexes. 
The analogous object for 2-complexes is an even 2-complex that is not a face-disjoint union of
two even 2-complexes. (That is, its faces cannot be partitioned into
two sets, each the face set of an even 2-complex.) We call such an object a {\bf circlet}. 

So a circlet is a minimal even 2-complex, in the sense that it contains no
proper even sub-2-complex.
{\it As a cycle is a minimal even 1-complex, a circlet is a minimal even 2-complex.}
Any connected surface is
a circlet. We'll soon see more exotic examples.

\begin{figure}[h]
\centering
\includegraphics{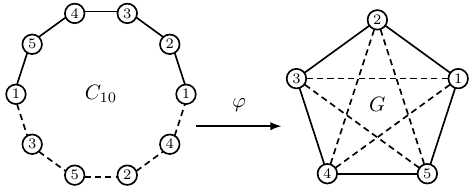}
\caption{An Euler tour $\varphi\colon C_{10}\to G$ in the graph $G= K_5$. Here
$V(G)=\{1,2, 3,4,5\}$, and each vertex $x$ of $C_{10}$ is labeled with $\varphi(x)$.
The solid and dashed paths in $C_{10}$ map to two edge-disjoint cycles in $G$, as illustrated.}
\label{Fig:K5Cover}
\end{figure}

\begin{figure}[t]
\centering
\includegraphics{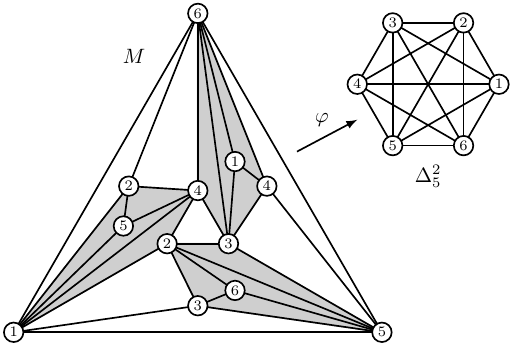}
\caption{An Euler cover $\varphi\colon M\to \Delta_5^{2}$ of the 2-skeleton of the
5-dimensional simplex, which has six vertices $1,2,3,4,5,6$, fifteen edges and twenty triangular faces.
Here $M$ is a triangulated sphere, with each vertex $x$ labeled
by $\varphi(x)$. The three shaded areas of $M$ map to three face-disjoint tetrahedra. The white areas of $M$ (including the unbounded face) map to an octahedron sharing no face with the three tetrahedra.}
\label{Fig:K6Cover}
\end{figure}

Let's now generalize the notion of an Euler tour.
An Euler tour in a graph $G$ can be described as a cellular map
$\varphi\colon C_n\to G$, where $C_n$ is a cycle and each edge of $G$ is the image of
exactly one edge of $C_n$.
See Figure~\ref{Fig:K5Cover}. Note that the domain $C_n$ is a connected
1-dimensional manifold. 
By strict analogy, we define an {\bf Euler cover} of a $2$-complex $K$ to be a cellular map
$\varphi\colon M\to K$, where $M$ is a connected 2-dimensional manifold (i.e., surface), and
each face of $K$ is the image (under $\varphi$) of exactly one face of $M$. So, whereas
an Euler tour charts the edges of a graph with a cycle, an Euler cover charts the faces 
of $K$ with a (connected) surface.
Figure~\ref{Fig:K6Cover} shows an Euler cover $\varphi: M \to K$, where $K$ is the
2-skeleton of the 5-dimensional simplex and $M$ is a triangulation of the sphere. 

We say a $2$-complex $K$ is {\bf strongly connected} if $K - V(K)$ is connected, that is,
if $K$ is connected and
removing its vertices does not disconnect it.
Here is our generalization of the Eulerian equivalence to 2-complexes.

\begin{theorem}
\label{Theorem:HK} For a strongly connected $2$-complex $K$, the following are equivalent:
\begin{itemize}[itemsep=0pt, align=right,topsep=4pt, labelwidth=1.5em,leftmargin=3.5em]
\item [{\em (i)}] $K$ is even,
\item [{\em (ii)}] $K$ is a face-disjoint union of circlets,
\item [{\em (iii)}] $K$ has an Euler cover.
\end{itemize}
\end{theorem}

Figure~\ref{Fig:K6Cover} illustrates the theorem. Here $K$ is the 2-complex consisting of the vertices, edges and
(triangular) faces of the 5-dimensional simplex.
The 2-complex $K$ is even because each edge belongs to 4 faces, and it is certainly strongly connected.
Figure~\ref{Fig:K6Cover} shows an Euler cover of $K$, and also illustrates how $K$ is a face-disjoint union
of circlets (three tetrahedral spheres and one octahedral sphere).

The next section gives examples of
circlets and proves that every circlet has an Euler cover. 
(And we'll see why every answer to Question 1 is correct.)
This is followed by a proof of Theorem~\ref{Theorem:HK},
and we conclude with examples and a discussion.

\section{Circlets and Their Euler Covers}
\label{Section:Circlets}

Recall that a circlet is a 2-dimensional generalization of a cycle, 
 an even 2-complex that is not a face-disjoint
union of two or more even 2-complexes. 

\begin{figure}[b]
\centering
\includegraphics[width=5in]{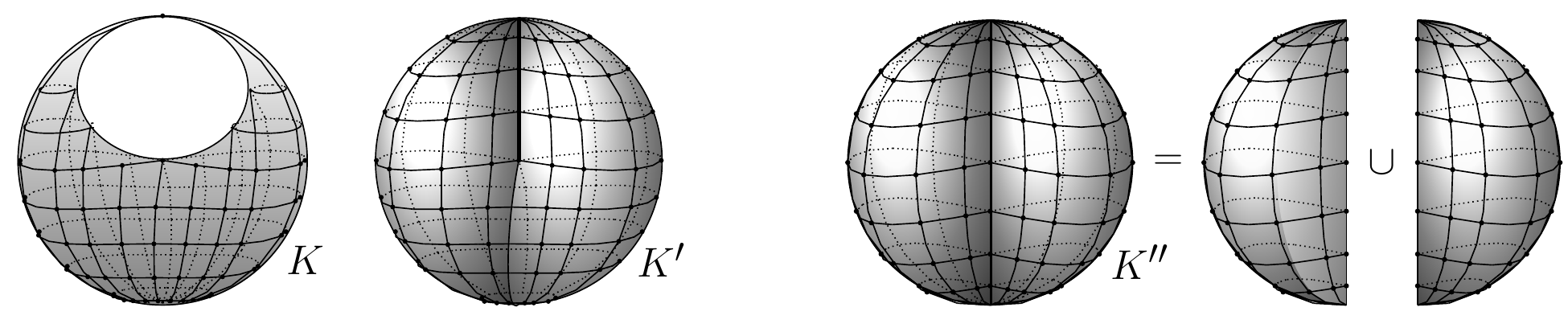}
\caption{Examples and non-examples of circlets. The pinched (or zipped) spheres $K$ and $K'$
are circlets, but the zipped sphere $K''$ is not.}
\label{Fig:CircEx}
\end{figure}

The most basic examples of circlets are the connected surfaces (orientable or non-orientable), as well as
connected surfaces with pinch points.
(See the pinched sphere $K$ in Figure~\ref{Fig:CircEx}.)
Also, ``zipping'' a connected surface along two paths can result in a circlet,
as is the case of $K'$ in Figure~\ref{Fig:CircEx}, in
which the portion of a great circle of a sphere is zipped on
itself. But zipping too far---as in $K''$
of Figure~\ref{Fig:CircEx}---can yield a 2-complex that is
not a circlet, as it is a face-disjoint union of two even complexes. 

The even complex $K$ of Figure~\ref{Fig:Klein} is another example of a circlet.
Indeed, suppose it were a face-disjoint union
$K=X\cup Y$ of two even 2-complexes. Say $X$ contains one of the faces having an edge
of degree~2. Then, to preserve evenness, $X$ must contain {\it all} faces having an edge of
degree~2. (This should be clear from the picture of $K$.) There is only one face remaining, namely the one
whose edges all have degree~4. But this must also belong to $X$, for otherwise $X$ is not even. As there are
no faces left over for $Y$, we admit that $K$ must be a circlet.

Our first proposition implies that every circlet is the face-disjoint image of a surface. (The converse is
false, as witnessed by $K''$ in Figure~\ref{Fig:CircEx}.)
We'll subsequently use
the proof's construction to show how all answers to
Question 1 are correct.

\begin{proposition}
Every circlet $K$ has an Euler cover $M\to K$. 
\label{Prop:CircletCover}
\end{proposition}

\begin{proof}
Our construction will fabricate $M$ by cutting $K$ along its edges, then re-gluing its faces so that two faces meet at each edge, as suggested in Figure~\ref{Fig:Explode}.

For each edge $e$ of $K$, let $\mathcal{F}_e$ be the set of faces of $K$ that have $e$ as an edge.
Select a partition $\mathcal{P}_e$ of  $\mathcal{F}_e$, each part of which
consists of exactly two faces (possible because $K$ is even).
Note that the choice of $\mathcal{P}_e$ is unique only if $\deg(e)=2$;
otherwise there are multiple ways to pick $\mathcal{P}_e$.
Given a face $f\in\mathcal{F}_e$, let $g_e(f)$ be the face paired with $f$
in $\mathcal{P}_e$. Thus $g_e^2(f)=f$,
so $g_e$ is an order-2 fixpoint-free permutation of $\mathcal{F}_e$.
We call $g_e$ a {\it gluing function} because below we will rip apart $K$ along its edges, then
glue any face $f$ to the face $g_e(f)$ along their common edge~$e$.

\begin{figure}[t]
\centering
\includegraphics{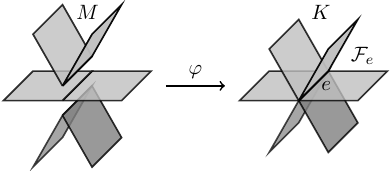}
\caption{The surface $M$ is obtained by unattaching the faces of $K$ and then re-gluing them along
their shared edges in pairs. If an edge $e$ of $K$ has degree $2k$, then it is the image of $k$ edges in $M$,
and the map $\varphi$ reattaches these $k$ edges to $e$.}
\label{Fig:Explode}
\end{figure}

Let $\widetilde{K}$ be the complex that is the disjoint union of all the faces of $K$.
For any face
$f$ of $\widetilde{K}$, there is an (injective) inclusion map $f\hookrightarrow K$. This induces an
``inclusion'' map $\iota: \widetilde{K}\to K$, which is not injective because if $x\in K$ is on the interior
of an edge~$e$, then its preimage $\iota^{-1}(x)$ has cardinality $\deg(e)$.

The gluing functions act on the faces of $\widetilde{K}$ in a way that mirrors their
effect on $K$. Specifically,
given a face $f\in\widetilde{K}$ and an edge $e$ of $f$, we agree that
 $g_e(f)$ is the face $f'\in\widetilde{K}$ for which $\iota(f')=g_{\iota(e)}(\iota(f))$.
 \begin{figure}[b]
\centering
\includegraphics{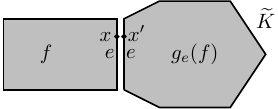}
\caption{Any face $f$ with edge $e$ is glued to $g_e(f)$ along their common edge.}
\label{FigTwoFaces}
\end{figure}

Now we make a surface by gluing the faces of $\widetilde{K}$
together edge-to-edge, two faces per edge. Identify any point $x$ on an 
edge $e$ of a face $f\in\widetilde{K}$ with its ``twin'' $x'$ on $g_e(f)$. That is,
such an $x\in e\subseteq f$ is identified with the $x'$ on $g_e(f)$ for which
$\iota(x)=\iota(x')$. See Figure~\ref{FigTwoFaces}. 
Let $q:\widetilde{K}\to M:=\widetilde{K}/\!\!\sim$ be the quotient map for
 this identification.
Then $M$ is a surface, and there is an induced cellular map $\varphi:M\to K$
making the following diagram commute. By construction, each face of $K$ is covered once.

\smallskip
{\centering
\begin{tikzpicture}[style=thick,scale=1]
\draw (0,0) node {$M$};
\draw (0,1) node {$\widetilde{K}$};
\draw (1.2,0) node {$K$};
\draw [->] (0,0.7)--(0,0.2);
\draw [->] (0.25,0)--(1,0);
\draw [->] (0.2,0.85)--(1,0.2);
\draw (-0.2,0.5) node {$q$};
\draw (0.6,-0.2) node {$\varphi$};
\draw (0.8,0.65) node {$\iota$};
\end{tikzpicture}
\par}

We claim that $M$ is connected. Say that $M$ is a disjoint
union $M=M_1\cup\cdots\cup M_n$ of connected surfaces. Then $K$ is a face-disjoint union
$K=\varphi(M_1)\cup\cdots\cup \varphi(M_n)$ of even 2-complexes. As $K$ is
a circlet, $n=1$, so $\varphi:M\to K$ is an Euler cover.
\end{proof}

To illustrate Proposition~\ref{Prop:CircletCover} (and its proof's construction) we next describe four
Euler covers of the circlet $K$ from Figure~\ref{Fig:Klein}. These four covers will show that each of the
four given answers to Question 1 is correct.

\begin{figure}[H]
\centering
\includegraphics{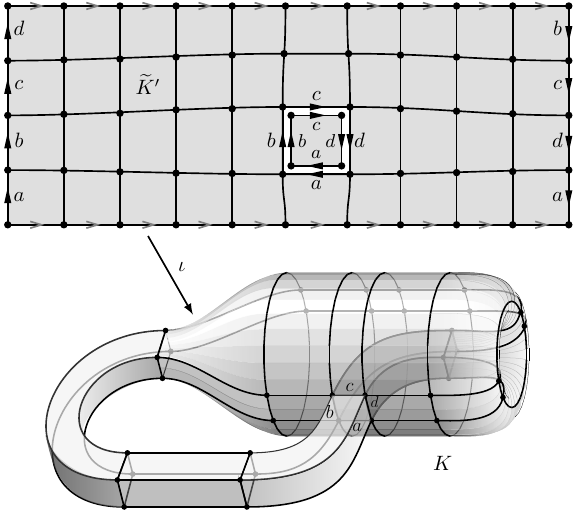}
\caption{The complexes $\widetilde{K}'$ and $K$.
We can think of $\widetilde{K}'$ as $K$ with all of its degree-4 edges 
unglued. Observe that $\widetilde{K}'$ is a tube (identify the upper and lower edges of the gray rectangle)
with a square cut out of it.
Euler covers of $K$ result from identifying 
the boundary edges of $\widetilde{K}'$ in pairs rather than in quadruples.}
\label{Fig:Kcut}
\end{figure}

The construction of an Euler cover $M\to K$ from Proposition~\ref{Prop:CircletCover}
begins with the disjoint union $\widetilde{K}$ of the faces of $K$, and the
selection of a gluing function $g_e$ for each edge of $K$. 
Notice that $K$ has four edges of degree~4 (labeled $a,b,c,d$ in the lower portion of
 Figure~\ref{Fig:Kcut}) and all other edges of $K$ have degree~2.
For an edge of $K$ having degree~4 (such as $a$)
there are $C(4,2)/2=3$ ways to choose the gluing function $g_a$. 
But for those edges $e$ of degree 2, there is only one choice for $g_e$.

So for any edge $e$ of $K$ having degree 2, 
the copies (in $\widetilde{K}$) of the two faces sharing~$e$ must be glued together in
$\widetilde{K}$ exactly as they are attached in $K$.
Denote by $\widetilde{K}'$ the complex $\widetilde{K}$ with all such face pairs
attached as they are attached in $K$.
It may be helpful to think of the edges of $K$ and $\widetilde{K}'$ as ``lines of glue'' that fuse
adjoining faces. Then $\widetilde{K}'$ is the space $K-abcd$, that is,
$K$ with the ``cycle of glue'' $abcd$ removed. (As faces are not yet glued along the edges
$a,b,c,d$.) Notice that $\widetilde{K}'$ is a tube with an isolated square cut
out of it, as shown at the top of Figure~\ref{Fig:Kcut}.
(To avoid fussy notation we use the symbols $a,b,c,d$ to stand for both edges of $K$
and their ``copies'' in  $\widetilde{K}'$.)

The boundary of $\widetilde{K}'$ comprises four cycles of length 4, each with edges labeled $a$, $b$, $c$ and~$d$. 
The complex $K$ results from gluing the four edges labeled $a$ together (with the indicated orientations), and likewise for the four labeled
$b$, $c$ and~$d$.
But the gluing functions $g_a,g_b,g_c,g_d$ encode instructions for gluing these edges
in pairs rather than quadruples. As there are 3 choices for each of these functions,
there are $3^4=81$ different gluings of the boundary edges of $\widetilde{K}'$,
each resulting in an Euler cover~$M$ of~$K$.
Below we construct four of them.

To begin, let's count the cells of $K$ and $M$. The grid form of $\widetilde{K}'$ makes this easy
for $K$. The complex $K$ has 
$4\cdot 10=40$ faces,
76 edges and 36 vertices.
(The four edges in $\widetilde{K}'$ labeled $a$ account for just one edge of $K$, etc.)
Say $M$ has $r$ faces, $q$ edges and $p$ vertices.
Then $r=40$ because an Euler cover $M\to K$ covers each of the 40 faces of $K$ once.
And $q=80$ because each of the 72 degree-2 edges of $K$ are covered once, but each of the
four degree-4 edges are covered {\it twice}.  A short 
argument\footnote{Outline: Each of the 32 vertices of $K$ that are not on the cycle $abcd$
must be covered exactly once by a vertex of $M$. So far, that accounts for 32 vertices
in $M$. We claim that each of the four vertices of the cycle $abcd$ is covered by either
{\it one or two} vertices of $M$
(so $36 \leq p \leq 40$).
Indeed, any vertex $x$ on $abcd$ touches eight faces of $K$,
so the preimage $\varphi^{-1}(x)$ is a set of
vertices of $M$ that collectively touch eight faces of $M$. If $|\varphi^{-1}(x)|>2$, then one vertex of
$\varphi^{-1}(x)$ touches at most $8/|\varphi^{-1}(x)|<3$ faces of $M$. Thus $M$ has a vertex
incident with only two faces of $M$, so these two faces share more than one edge.
Then $\varphi$ maps
these faces to two faces of $K$ that share more than one edge. But $K$ has no such pair of faces, so
$|\varphi^{-1}(x)|\leq 2$, proving the claim.}
shows  $36\leq p\leq 40$. 
So $M$ has Euler characteristic
\begin{equation}
\chi(M)\;=\;p-q+r \;=\; p-40\;\in \;\{0,-1, -2, -3,-4\}.
\label{Eqn:EC}
\end{equation}

This tags $M$ as one of eight surfaces:
If $\chi(M)$ is $0$, then $M$ is a torus or Klein bottle.
If $\chi(M)$ is $-1$ or $-3$, then $M$ is a sphere with three or five cross-caps, respectively.\linebreak
If $\chi(M)$ is $-2$, then $M$ is a double torus or a sphere with four cross-caps.
Finally,\linebreak
if $\chi(M)$ is $-4$, then $M$ is either a triple torus or a sphere with six cross-caps.

\begin{figure}[H]
\centering
\includegraphics{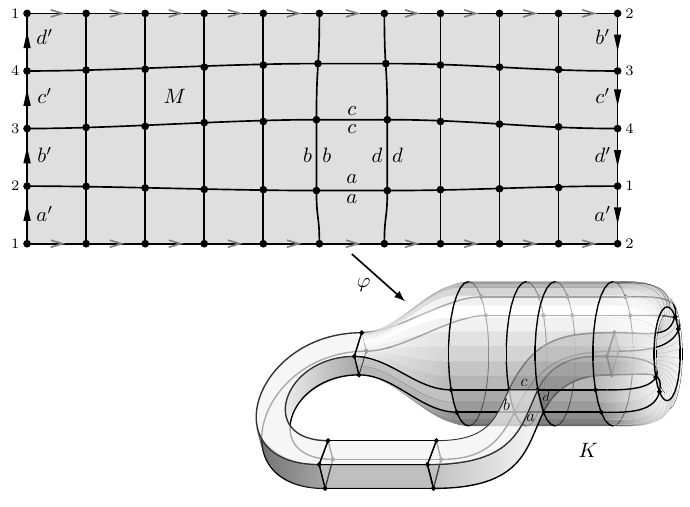}
\caption{{\bf Euler cover of K by a Klein bottle:} Make a surface $M$ by gluing the isolated square of $\widetilde{K'}$ 
(Figure~\ref{Fig:Kcut}) to the shaft of the
tube, as shown here. Then glue the left end of the tube to the right end, so
that $a'$ is glued to $a'$, $b'$ to $b'$, etc.
The cycle $a'b'c'd'$ is now distinct from the cycle $abcd$, and
the vertices 1,2,3,4 of the cycle $a'b'c'd'$ are distinct from the vertices of $abcd$.
The resulting surface $M$ is a
Klein bottle.
The Euler cover $\varphi:M\to K$ maps $a',a\in E(M)$ to $a\in E(K)$.
Also $\varphi(b')=\varphi(b)=b \in E(K)$, etc.}
\label{Fig:Kcover}
\end{figure}

\begin{figure} [H]
\includegraphics{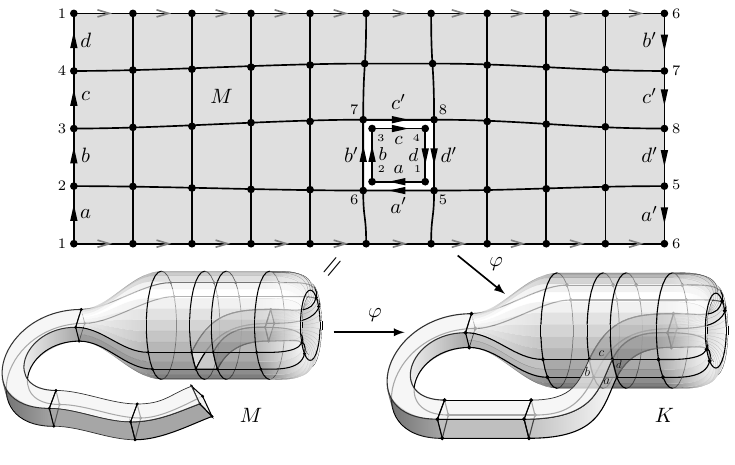}
\caption{{\bf Euler cover of K by a torus:} 
Glue the tube's left boundary cycle $abcd$ to the isolated square,
thus capping the left end of the tube. Glue the tube's right boundary $a'b'c'd'$
to the boundary of the square hole, as shown. This creates a handle,
resulting in a torus $M$.
Informally, $M$ is $K$ with its neck detached and capped with the square. So $M$ is a sphere with a tunnel through it, that is, a torus. The Euler cover $M\to K$ re-attaches the neck.}
\label{Fig:T1cover}
\end{figure}

\begin{figure}[H]
\includegraphics{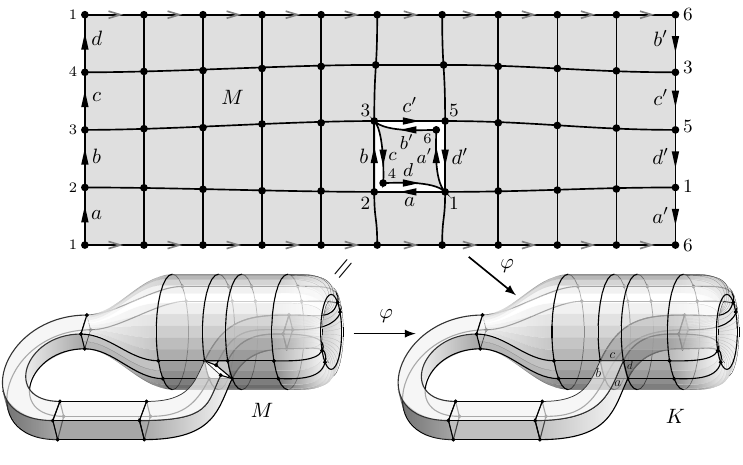}
\caption{{\bf Euler cover of K by a double torus:} Form $M$ by attaching the isolated square to the
shaft of the tube at two diagonal vertices. Flip the square along the diagonal.
Glue
the tube's left boundary $abcd$ to the hole bounded by $abcd$ (thus creating a handle).
Glue the
right boundary $a'b'c'd'$  to  the hole bounded by $a'b'c'd'$ (a second handle).
This is also illustrated on the bottom. In~$M$, the isolated square pivots to vent
a tunnel that exits to the right. 
(On the other side of the neck it vents a hole that connects the interior of the neck to the inner
chamber of the bottle.)
The Euler cover $\varphi$ closes the vent.
}
\label{Fig:T2Cover}
\end{figure}

\begin{figure}[h]
\includegraphics{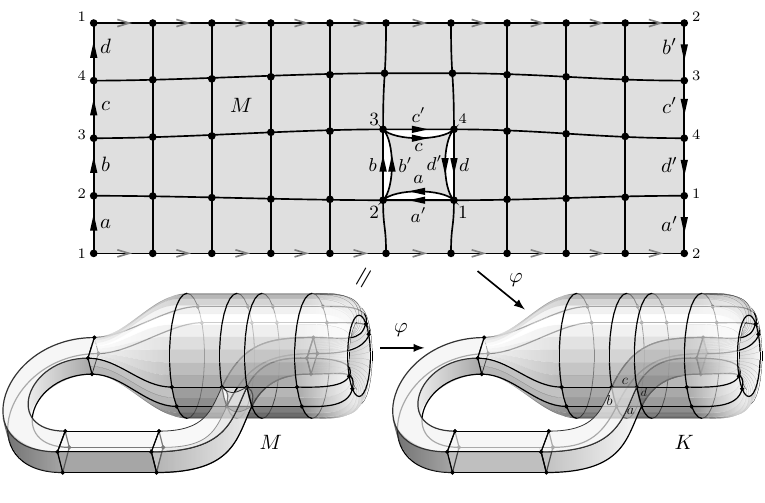}
\caption{{\bf Euler cover of K by a triple torus:} 
Attach the isolated square from Figure~\ref{Fig:Kcut} to the tube at its four corners and glue the arcs 
in pairs as shown. The resulting surface $M$ has 36 vertices, the same number of vertices as $K$.
(All vertices of $M$ have degree~4, except those
labeled $1,2,3,4$, each of which has degree~8.) Equation~(\ref{Eqn:EC}) gives $\chi(M)=-4$,
so $M$ can only be a triple torus or a sphere with six cross-caps. 
Inspection reveals no M\"obius strips, so $M$ is a triple torus.
This may be more easily visualized as  drawn at the bottom.
Here the isolated square is saddle-shaped. Where the neck meets the main bottle cylinder there are
two crescent-shaped holes (bounded by dicycles $aa'$ and $cc'$) that give passage into the main tunnel that
exits to the right. (And the two dicycles $bb'$ and $dd'$ are two openings through which the
interior of the neck passes into the interior of the main cylinder.)
}
\label{Fig:T3Cover}
\end{figure}

Each of these possibilities is attainable with a judicious choice of
gluing functions $g_a$, $g_b$, $g_c$ and $g_d$.
Figures~\ref{Fig:Kcover}--\ref{Fig:T3Cover}
describe four of them, the answers to Question 1. We leave the remaining four cases as exercises.

\pagebreak

\vspace*{-0.1in}

\pagebreak
 
\section{Proof of the Main Theorem}
\label{Section:EulerExist}

We can now prove Theorem~\ref{Theorem:HK}. Suppose $K$ is a strongly connected 2-complex.
Below we prove (iii)$\Rightarrow$(i)$\Rightarrow$(ii)$\Rightarrow$(iii).
(Strong connectivity is required only for (ii)$\Rightarrow$(iii). 
In fact, the reader can easily verify that (iii)$\Rightarrow$(i)$\Leftrightarrow$(ii) hold without it.)

\begin{proof}
{\bf (iii)$\Rightarrow$(i)} (If $K$ has an Euler cover, then $K$ is even.)

\noindent
If a 2-complex $K$ has an Euler cover $M\to K$, then any edge $e$ of $K$ has degree $2k$, where $e$
is the image of $k$ edges in $M$.

\medskip
\noindent
{\bf (i)$\Rightarrow$(ii)} (If $K$ is even then $K$ a face-disjoint union of circlets.)

\noindent
Suppose $K$ is even. If $K$ is not the face-disjoint union of two even 2-complexes, then $K$ is
a circlet and we are done. Otherwise  $K=X\cup Y$ is a face disjoint union of two even 2-complexes
$X$ and $Y$. If neither $X$ nor $Y$ are face-disjoint unions of even 2-complexes, then they
are circlets, and we are done. Otherwise at least one of them decomposes into a face-disjoint union of
even 2-complexes. In a finite number of steps this process decomposes $K$  into a face-disjoint union of circlets.

\pagebreak
\noindent
{\bf (ii)$\Rightarrow$(iii)} (If $K$ is a face-disjoint union of circlets, then $K$ has an Euler cover.)

\noindent
Suppose $K$ is a face-disjoint union  $K=K_1\cup \cdots \cup K_n$ of circlets.
Proposition~\ref{Prop:CircletCover} guarantees an Euler cover
$\varphi_i:M_i\to K_i$ of each circlet. This induces a cellular map $\varphi:(M_1\cup \cdots \cup M_n)\to K$
having all the properties of an Euler cover except that the domain
surface is disconnected if $n>1$.
In what follows we explain how to splice the $M_i$ together to get an Euler cover
$M\to K$.

\begin{figure}[b]
\centering
\includegraphics{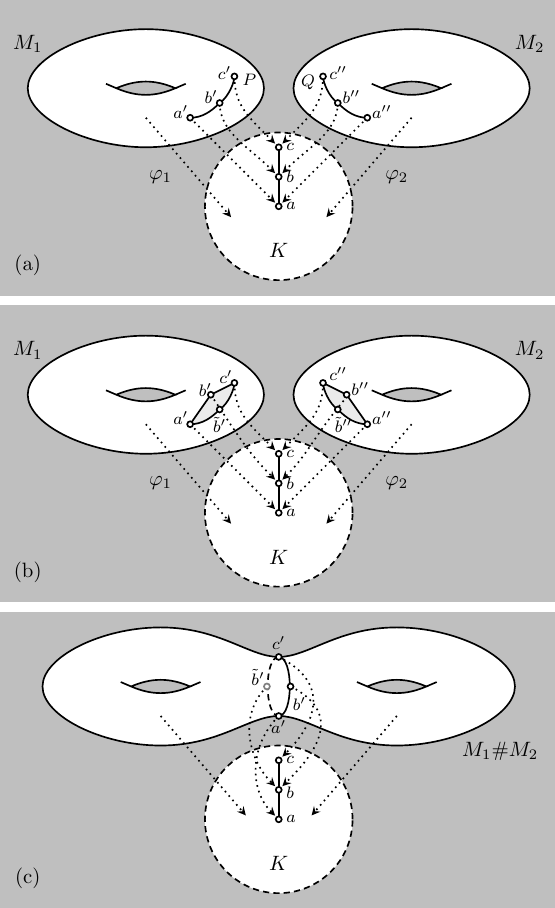}
\caption{The construction for splicing covers.}
\label{FIG:EulerCover2}
\end{figure}

Now, $K$ must have an edge $ab$ that is the image of an edge $a'b'\in E(M_1)$ and 
also the image of an edge $a''b''$ of some other $M_j$.
Indeed, if this were not so, then no edge
of $\varphi(M_1)$ would be an edge of any $\varphi(M_j)$ and hence $\varphi(M_1)$ could be
disconnected from the rest of $K$ by removing all vertices of $K$, violating strong connectivity
of $K$. 
There will be no harm in assuming  $M_j=M_2$. See Figure~\ref{FIG:EulerCover2}.
Choose the labeling so that $\varphi_1(a')=a=\varphi_2(a'')$ and $\varphi_1(b')=b=\varphi_2(b'')$.

In general, $M_1$ and $M_2$ will have paths $P$ (with initial edge $a'b'$) and $Q$ (with initial edge  $a''b''$),
respectively, of common length at least~1, for which
$\varphi_1(P)=\varphi_2(Q)$. (See Figure~\ref{FIG:EulerCover2}.)
Slit $M_1$ along~$P$, opening a hole, making~$M_1$ a surface with boundary.
In a like manner slit~$M_2$ along $Q$, opening a hole.

Now glue $M_1$ and $M_2$ together along their boundaries, in such a way that
any boundary point $x'$ of $M_1$ is glued to a boundary point $x''$ of $M_2$ for which
$\varphi_1(x')=\varphi_2(x'')$, as indicated in Figure~\ref{FIG:EulerCover2}(c).
This yields a new surface $M_1\# M_2$, a connected sum of $M_1$ and $M_2$.
Because $\varphi_1$ and $\varphi_2$ agree on the identified points, they induce a cellular map
$M_1\# M_2\to K$ that is an Euler cover of its range.
(If $P$ and $Q$ have length 1, then this construction introduces a pair of parallel edges in $M_1\# M_2$.)

Performing this operation reduces the number of components of 
$M_1\cup\cdots \cup M_n$ by one, while preserving the salient features of an Euler cover.
We can perform this operation in a sequence, reducing the number of components of our 
domain by one in each iteration. We eventually  arrive at just one component $M$ and we are done.
\end{proof}

\section{Decomposition into Circlets}
\begin{figure}[b]
\centering
\includegraphics[width=2.2in]{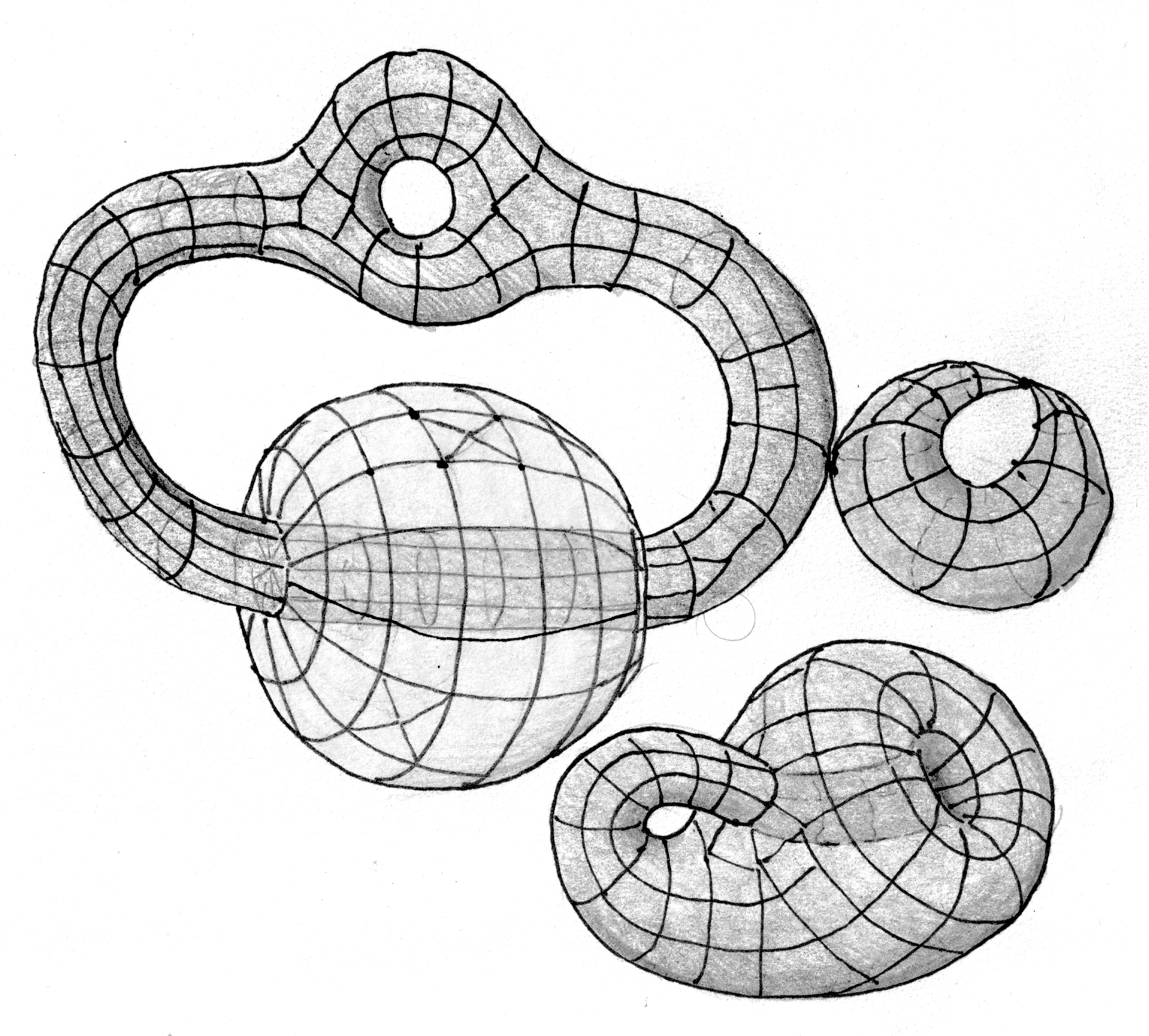}
\caption{Every even complex is a face-disjoint union of circlets.}
\label{Fig:Surfaces}
\end{figure}

We have hitherto concentrated mainly on Euler covers of even 2-complexes, that is, the equivalence (i)$\Leftrightarrow$(iii)
of Theorem~\ref{Theorem:HK}. Let's now turn to the equivalence (i)$\Leftrightarrow$(ii),
guaranteeing that any even 2-complex is a face-disjoint union of circlets.
This is illustrated (somewhat whimsically) in Figure~\ref{Fig:Surfaces}, showing an even 2-complex that
is a face-disjoint union of a double torus, a sphere, an image of Klein bottle and a pinched sphere.
This decomposition is not unique; this 2-complex also decomposes as two tori,
an image of a Klein bottle and a pinched sphere.

In~\cite{PLATONIC} we describe explicit circlet decompositions of all even 2-skeletons of the $n$-dimensional
platonic polyhedra. For example, consider the 2-skeleton of an odd-dimensional simplex.

The {\bf n-dimensional simplex}  $\Delta_n$ is the convex hull
of the standard basis elements $\{{\bf e}_1, {\bf e}_2, \ldots, {\bf e}_{n+1}\}\subseteq\mathbb{R}^{n+1}$. 
The $2$-cells (or faces) of $\Delta_n$ are precisely the convex hulls of the triples
$\{{\bf e}_{i_1}, {\bf e}_{i_2}, {\bf e}_{i_{2}}\}$ for $\{i_1,i_2,i_3\}\subseteq \{1,2,\ldots, n{+}1\}$.
The 2-skeleton $\Delta_n^2$ of $\Delta_n$ is
the 2-complex consisting of the $n{+}1$ vertices,  $\binom{n+1}{2}$ edges, and $\binom{n+1}{3}$ faces
of $\Delta_n$. It is even precisely when $n$ is odd, as each edge belongs to
$n{-}1$ faces.

\begin{figure}[h]
\centering
\includegraphics{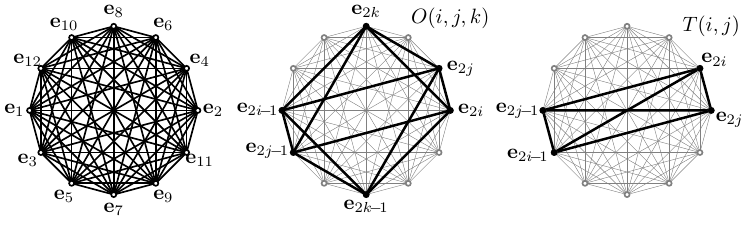}
\caption{Decomposition of $\Delta_n^2$ into circlets. Left: $\Delta_n$. Middle: $O(i,j,k)$. Right $T(i,j)$.}
\label{FigOijkTij}
\end{figure}

Here is a recipe for decomposing $\Delta_n^2$ into spherical circlets when $n$ is odd. 
Draw $\Delta_n$ as a regular $(n{+}1)$-gon with an edge connecting each vertex pair.
Label the vertices in a counterclockwise
sequence ${\bf e}_2$,  ${\bf e}_4$,  ${\bf e}_6$, $\ldots$  ${\bf e}_{n+1}$, followed by
${\bf e}_1$,  ${\bf e}_3$,  ${\bf e}_5$, $\ldots$  ${\bf e}_{n}$, so that ${\bf e}_{2i}$ is always
the antipode (relative to the drawing) of ${\bf e}_{2i-1}$. See Figure~\ref{FigOijkTij}, left.

For each triple $\{2i, 2j, 2k\}\subseteq \{2,4,6,\ldots, n{+}1\}$, 
$\Delta_n$ has eight triangle faces whose vertices are the eight triples from
$\{{\bf e}_{2i}, {\bf e}_{2i-1}\}\times\{ {\bf e}_{2j}, {\bf e}_{2j-1}\}\times\{ {\bf e}_{2k}, {\bf e}_{2k-1}\}$.
These eight triangles are shown in the middle
drawing in Figure~\ref{FigOijkTij}. They constitute the faces of
a 2-complex that we shall denote as $O(i,j,k)$. Observe that $O(i,j,k)$
is homeomorphic the boundary of an octahedron.

Also, for $\{2i, 2j\}\subseteq \{2,4,6,\ldots, n{+}1\}$,  the convex
hull of $\{{\bf e}_{2i},  {\bf e}_{2i-1}, {\bf e}_{2j}, {\bf e}_{2j-1}\}$
is a tetrahedron whose boundary $T(i,j)$ is six triangular faces
of $\Delta_n$ (Figure~\ref{FigOijkTij}, right).

Thus in $\Delta_n^2$ we have $\binom{(n+1)/2}{3}$ octahedral spheres $O(i,j,k)$ and $\binom{(n+1)/2}{2}$
tetrahedral spheres $T(i,j)$. Each face of $\Delta_n$ belongs to precisely one of these spheres,
so we have realized $\Delta_n^2$ as a face-disjoint union of circlets, all of which are spheres. 

We can create an Euler cover of $\Delta_n^2$ with the construction in the proof of our main theorem. 
For example, see the Euler cover and decomposition of $\Delta_5^2$ in Figure~\ref{Fig:K6Cover}.

In~\cite{PLATONIC} we also show that the 2-skeleton of any cross-polytope is a face-disjoint union of
octahedral spheres, and that odd-dimensional hypercubes are face-disjoint unions of tori and spheres.
In~\cite{GENUS}, the 2-skeleton of the $n$-dimensional hypercube (for odd $n$) is decomposed into 
$(n{-}1)/2$ isometric multi-holed tori, into which the 1-skeleton (i.e., the hypercube graph) is
simultaneously minimally 2-cell embedded.

An interesting question is whether the 2-skeleton of every odd-dimensional hypercube decomposes into
spheres. In~\cite{SPHERES} we use design theory to show that this is possible in dimensions 
$n\equiv 1 \mbox{ or } 3$ (mod 6), as well as in other sporadic dimensions. Whether this is
possible in every dimension is unknown.

\section{Conclusion}
We've shown that the three conditions of Euler's theorem (Theorem~\ref{Theorem:Euler})
for graphs also make sense in two dimensions and remain equivalent.  Using the one-dimensional case as a model, we were led to generalizations of cycles and Euler tours.  
While an Euler tour can be seen as a
cycle (connected 1-manifold) that intersects itself at points,
our notion of Euler cover of a 2-complex has
a 2-manifold intersecting itself only in its
1-skeleton.

Although this article focuses on 2-complexes, our results do generalize to $n$-dimensions, with some caveats. 
Certainly the notions of evenness and circlets generalize nicely: an $n$-complex is {\bf even} if each
$(n{-}1)$-cell belongs to a positive even number of $n$-cells; and an {\bf n-circlet} is an even $n$-complex
that is not a $n$-cell-disjoint union of two even $n$-complexes. With these definitions,
an $n$-complex is even if and only if it decomposes into a $n$-cell-disjoint union of $n$-circlets.

But the notion of an $n$-dimensional Euler tour is more subtle. 
The constructions presented here for $n=2$ do not work so well if $n>2$. The heart of the issue
stems from the fact that when polygonal faces are glued together edge-to-edge, faces arrange themselves 
cyclically around vertices, so that each vertex has a neighborhood homeomorphic to a disk (thus
yielding a surface, or 2-dimensional manifold). For $n>2$ there is no such cyclic arrangement, and gluing $n$-cells together along
$(n{-}1)$-cells may not result in a $n$-dimensional manifold, but rather a pseudo-manifold. 
For a detailed treatment, see \cite{COMPLEX}.

\begin{acknowledgment}
We thank the referees and the editor. RHH is supported by Simons Foundation Collaboration Grant for Mathematicians 523748.
\end{acknowledgment}

\summary{A famous theorem in graph theory---originating with Euler---characterizes connected even-degree graphs
as (1) those graphs that admit an Euler tour, and (2) those connected graphs that decompose as a face-disjoint union of cycles. 
We explore a 2-dimensional generalization of this theorem, with graphs (i.e., 1-complexes) replaced by
2-complexes. This entails an interesting generalization of cycles, and the introduction
of the notion of a ``2-dimensional Euler tour.''}

\begin{biog}
\item[Richard Hammack] (MR Author ID: 649820)
is a professor of mathematics at Virginia Commonwealth University,
working mostly in graph theory and combinatorics. He is the author of
{\it Book of Proof}, an open proofs textbook, and coauthor of
{\it Handbook of Product Graphs} (with W.\ Imrich and S.\ Klav\v{z}ar).
\end{biog}

\begin{biog}
\item[Paul Kainen] (MR Author ID: 96855)
is adjunct professor of mathematics at Georgetown University. He works mostly in geometric graph theory, category theory and neural networks,
and he is co-author of {\em The Four Color Theorem: Assaults and Conquest} (with T. L. Saaty).
\end{biog}

\end{document}